# SOME REMARKS ON BOUNDARY OPERATORS OF BESSEL EXTENSIONS

JESSE GOODMAN AND DANIEL SPECTOR

ABSTRACT. In this paper we study some boundary operators of a class of Bessel-type Littlewood-Paley extensions whose prototype is

$$\Delta_x u(x,y) + \frac{1-2s}{y}\frac{\partial u}{\partial y}(x,y) + \frac{\partial^2 u}{\partial y^2}(x,y) = 0 \qquad \text{for } x \in \mathbb{R}^d, y > 0,$$
$$u(x,0) = f(x) \qquad \text{for } x \in \mathbb{R}^d.$$

In particular, we show that with a logarithmic scaling one can capture the failure of analyticity of these extensions in the limiting cases $s = k \in \mathbb{N}$.

## 1. Introduction and Main Results

In the paper [1], Caffarelli and Silvestre showed that for $s \in (0,1)$, the fractional Laplacian

$$(1.1) \qquad (-\Delta)^s f(x) = c_{d,s} \text{ p.v.} \int_{\mathbb{R}^d} \frac{f(x)-f(y)}{|x-y|^{d+2s}} \, dy$$

can be obtained via a Dirichlet to Neumann-type operator of the elliptic partial differential equation

$$(1.2) \qquad \Delta_x u(x,y) + \frac{1-2s}{y}\frac{\partial u}{\partial y}(x,y) + \frac{\partial^2 u}{\partial y^2}(x,y) = 0 \qquad \text{for } x \in \mathbb{R}^d, y > 0,$$
$$u(x,0) = f(x) \quad \text{for } x \in \mathbb{R}^d.$$

Namely, the solution $u$ of (1.2) satisfies the boundary asymptotics

$$(1.3) \qquad \lim_{y \downarrow 0} \frac{u(x,y) - u(x,0)}{y^{2s}/2s} = -\frac{\Gamma(1-s)}{2^{2s-1}\Gamma(s)}(-\Delta)^s f(x).$$

Caffarelli and Silvestre's establishment of such a formula makes an important contribution in connecting ideas that had been explored in more classical papers concerning the problem (1.2) but were not fully realized. Indeed, Molčanov and Ostrovskiĭ had in [7, Theorem 2] given a description of the Dirichlet to Neumann-type map for the problem (1.2) in a probabilistic language by finding the trace process of the degenerate diffusion (1.2), while Marias [6, p. 320, equation (3.1)] had derived the Poisson-Bessel kernel for (1.2),

$$(1.4) \qquad u(x,y) = \frac{\Gamma(\frac{d}{2}+s)}{\pi^{d/2}\Gamma(s)} \, y^{2s} \int_{\mathbb{R}^d} \frac{f(x-z)}{(|z|^2+y^2)^{d/2+s}} \, dz,$$

in his paper concerning Littlewood-Paley-Stein theory for the equation (1.2).

Given the well-posedness of the PDE (1.2) for all $s \in \mathbb{R}^+$, a natural question is whether one can obtain a boundary operator in the spirit of (1.3) for larger values of $s$. In fact for $s \in (0, \frac{d}{2}) \setminus \mathbb{N}$ such a result has been obtained by Chang and Gonzalez [2] (see also Yang [12]), while Roncal and Stinga [9] have extended their result to the case $s \in \mathbb{R}^+ \setminus \mathbb{N}$.

---

In this paper we undertake the problem of obtaining the boundary behavior for the PDE (1.2) in this spirit when $s$ is a non-negative integer, thus completing the picture for any $s \geq 0$. We begin with $s = 0$. In this case, the boundary-value problem (1.2) is not well-posed. We can nevertheless obtain a solution $u(x, y)$ parametrized by a function $f(x)$ on the boundary by formally replacing $\Gamma(s)$ with 1 in the Poisson-Bessel formula (1.4) for $s = 0$.

**Theorem 1.1.** *Let $f \in C_c^\infty(\mathbb{R}^d)$ and define the function $u$ by*

$$u(x, y) = \frac{\Gamma(d/2)}{\pi^{d/2}} \int_{\mathbb{R}^d} \frac{f(x - z)}{(|z|^2 + y^2)^{d/2}} \, dy \qquad x \in \mathbb{R}^d, y > 0.$$

*Then $u$ is a solution of*

$$\Delta_x u(x, y) + \frac{1}{y}\frac{\partial u}{\partial y}(x, y) + \frac{\partial^2 u}{\partial y^2}(x, y) = 0 \qquad x \in \mathbb{R}^d, y > 0$$

*and has the boundary asymptotics*

$$\lim_{y \downarrow 0} \frac{1}{\ln(y)} u(x, y) = -2f(x).$$

In general, $u(x, y)$ fails to be analytic as a function of $y$ near $y = 0$. Thus Theorem 1.1 shows that in the case $s = 0$, the leading-order behavior of the solution $u$ is singular, proportional to $\ln(y)$, whereas (1.3) shows that, for $0 < s < 1$, $u(x, y)$ equals to leading order the term $u(x, 0)$, followed by a singular term proportional to $y^{2s}$ (except in the case $s = 1/2$, when this term is itself analytic in $y$). For larger $s$, we can expect the asymptotic expansion of $u(x, y)$ to contain many analytic terms of lower order than the first singular term.

To see this concretely, note that we can represent the solution $u$ formally as

(1.5) $$u(x, y) = \int_{\mathbb{R}^d} \hat{f}(\xi) e^{2\pi i \xi \cdot x} \varphi(2\pi|\xi|y) \, d\xi$$

where

$$\varphi(t) := \frac{1}{\Gamma(s)} \int_0^\infty v^{-s} \exp\left(-\frac{1}{v} - \tfrac{1}{4}t^2 v\right) \frac{dv}{v}.$$

Notice this function $\varphi$ is the explicit solution of the ODE considered in Caffarelli and Silvestre's paper on p. 1250. Expand the exponential function in the definition of $\varphi(2\pi|\xi|y)$ as a power series in $y$,

$$\exp(-\tfrac{1}{4}(2\pi|\xi|y)^2 v) = \sum_{n=0}^\infty \frac{(2\pi^2|\xi|^2 y^2)^n (\frac{-v}{4})^n}{n!}$$

and define

$$c_n := \frac{1}{\Gamma(s)} \int_0^\infty \frac{1}{n!}(\frac{-v}{4})^n v^{-s} \exp\left(-\frac{1}{v}\right) \frac{dv}{v} = \frac{(-1)^n \Gamma(s - n)}{2^{2n} n! \Gamma(s)}.$$

Then $c_n$ is finite whenever $n < s$, and we may perform the integration over $v$ termwise for these values. Thus we expect the asymptotic expansion of $u(x, y)$ to begin with $\lceil s \rceil$ analytic terms before reaching the first singular term:

$$u(x, y) = f(x) + y^2 c_1(-\Delta) f(x) + (y^2)^2 c_2(-\Delta)^2 f(x)$$
$$+ \cdots + (y^2)^{\lceil s \rceil - 1} c_{\lceil s \rceil - 1}(-\Delta)^{\lceil s \rceil - 1} f(x)$$
$$+ \int_{\mathbb{R}^d} \hat{f}(\xi) e^{2\pi i \xi \cdot x} \frac{1}{\Gamma(s)} \int_0^\infty \sum_{n=\lceil s \rceil}^\infty \frac{(2\pi^2|\xi|^2 y^2)^n (\frac{-v}{4})^n}{n!} v^{-s} \exp\left(-\frac{1}{v}\right) \frac{dv}{v}.$$

Now, as in (1.3) and Theorem 1.1, we are interested in the failure of analyticity which takes place in the last term. Therefore we introduce a differentiation operator that allows us to remove the lower order terms in $y^2$, which is

$$\frac{\partial}{\partial(y^2)} = \frac{1}{2y}\frac{\partial}{\partial y}.$$

Notice this is none other than the operator that arises in the induction argument of Chang and Gonzalez in [2], namely

$$(1.6) \qquad Tu(x,y) := \frac{1}{y}\frac{\partial u}{\partial y}(x,y),$$

except for a factor of two. When $s$ is not an integer, say $s = k + \sigma$ where $k = \lfloor s \rfloor \in \mathbb{N} \cup \{0\}$ and $\sigma \in (0,1)$, one can apply this operator $k$ times to remove most of the lower-order terms, whereupon an argument similar to the computation (1.3) enables one to obtain

$$(1.7) \qquad \lim_{y\downarrow 0} \frac{T^k u(x,y) - T^k u(x,0)}{y^{2\sigma}/2\sigma} = \frac{(-1)^{k-1}\Gamma(1-\sigma)}{2^k \Gamma(s) 2^{2\sigma-1}}(-\Delta)^s f(x),$$

which is one of the results obtained in the papers [2, 9, 12].

This analysis motivates our next result, which concerns the case $s = k \in \mathbb{N}$. In this regime we have the following theorem which completes the picture concerning the equation (1.2) for $s > 0$.

**Theorem 1.2.** *Let $s = k \in \mathbb{N}$ and suppose $u$ solves (1.2) with boundary data $f \in C_c^\infty(\mathbb{R}^d)$. Then*

$$\lim_{y\downarrow 0} \frac{T^k u(x,y)}{\ln(y)} = \frac{(-1)^{k-1}}{2^{k-1}\Gamma(k)}(-\Delta)^k f(x).$$

Our proofs of Theorems 1.1 and 1.2 are via explicit computation, as it allows us to connect some of the constants introduced in the various papers [1, 2, 11] (see Section 2). However, as one observes from the paper of Stinga and Torrea [11], not only does the constant arising in the computation (1.3) not depend on the dimension, it does not even depend on the fact that one performs an extension of the Laplacian. More precisely, consider the PDE

$$(1.8) \qquad \begin{aligned} -Lu(x,y) + \frac{1-2s}{y}u_y(x,y) + u_{yy}(x,y) = 0 \text{ in } \mathbb{R}_+^{d+1} \\ u(x,0) = f(x) \text{ on } \mathbb{R}^d, \end{aligned}$$

for $L$ a non-negative elliptic operator, densely defined, and self-adjoint on $L^2(\mathbb{R}^d)$. Then for $s \in (0,1)$ Stinga and Torrea [11] have shown one has the analogue to (1.3), i.e.

$$\lim_{y\downarrow 0} \frac{u(x,y) - u(x,0)}{y^{2s}/2s} = -\frac{\Gamma(1-s)}{2^{2s-1}\Gamma(s)} L^s f(x),$$

while for $s \in \mathbb{R}^+ \setminus \mathbb{N}$ Roncal and Stinga [9] have established boundary asymptotics in the spirit of (1.7).

In this paper we take up the analysis of $s = k \in \mathbb{N} \cup \{0\}$ for more general operators. To motivate our approach, let us first observe that an alternative formulation of the function $u(x,y)$ from Theorem 1.1 can be made in terms of the heat semigroup

$$e^{t\Delta}f(x) = \int_{\mathbb{R}^d} \frac{1}{(4\pi t)^{d/2}} e^{-|z-x|^2/4t} f(z)\,dz.$$

Precisely, we have

**Proposition 1.3.** *Let $f \in C_c^\infty(\mathbb{R}^d)$. Then*
$$\frac{\Gamma(d/2)}{\pi^{d/2}} \int_{\mathbb{R}^d} \frac{f(z)}{(|x-z|^2+y^2)^{d/2}} \, dz = \int_0^\infty e^{t\Delta} f(x) e^{-\frac{y^2}{4t}} \frac{dt}{t}.$$

We note in particular that the relation between the function $u(x,y)$ and the semigroup $e^{t\Delta}f(x)$ does not depend on the dimension $d$, while such a representation allows us to cast the convergence from Theorem 1.1 as the more abstract formula

$$\lim_{y \downarrow 0} \frac{1}{\ln(y)} \int_0^\infty e^{t\Delta} f(x) e^{-\frac{y^2}{4t}} \frac{dt}{t} = -2f(x).$$

This perspective suggests that a good framework for posing our question of boundary operators for more general equations than the Laplacian is that of semigroups. Thus in our final two theorems, we consider a strongly continuous semigroup of bounded linear operators $S(t)$ from a Banach space $\mathcal{F}$ to itself. We write $\|\cdot\|$ for the norm associated to $\mathcal{F}$. Let $-L$ denote the infinitesimal generator of $S(t)$ and let $\mathcal{D}(L) \subset \mathcal{F}$ denote its domain. That is,

$$(1.9) \qquad -Lf = \lim_{t \downarrow 0} \frac{S(t)f - f}{t}$$

and $\mathcal{D}(L)$ is the set of $f \in \mathcal{F}$ for which the limit in (1.9) exists as an element of $\mathcal{F}$. We assume that the semigroup is bounded in the sense that

$$(1.10) \qquad \sup_{t \geq 0} \|S(t)f\| < \infty \qquad \text{for each } f \in \mathcal{F}.$$

We will use the semigroup $S(t)f$ to construct functions $u(f, y)$ that satisfy

$$(1.11) \qquad -Lu(f,y) + \frac{1-2s}{y} \frac{\partial u}{\partial y} + \frac{\partial^2 u}{\partial y^2} = 0, \quad y > 0.$$

When $\mathcal{F}$ is a function space and $f = f(x)$ is a function on $\mathbb{R}^d$, $u(f, y)$ is the analogue of $u(\cdot, y)$ in Theorems 1.1–1.2.

The assumptions (1.9)–(1.10) cover a large class of the operators included in the work of Stinga and Torrea [11]. For example, one can take

$$-L = \sum_{i=1}^d \sum_{j=1}^d \frac{\partial}{\partial x_i}\left(a_{ij}(x) \frac{\partial}{\partial x_j}\right),$$

or more generally

$$-L = \sum_{i=1}^d \sum_{j=1}^d a_{ij}(x) \frac{\partial^2}{\partial x_i \partial x_j} + \sum_{i=1}^d b_i(x) \frac{\partial}{\partial x_i},$$

where the coefficients $a_{ij}(x)$ and $b_i(x)$ are smooth bounded functions of $x \in \mathbb{R}^d$ with $a_{ij}$ having bounded first and second derivatives, $b_i$ with bounded first derivatives, and the matrix with entries $a_{ij}(x)$ is symmetric and positive definite for each $x \in \mathbb{R}^d$. Then for $f \in L^2(\mathbb{R}^d)$, there exists a unique solution of the parabolic equation

$$\frac{\partial g}{\partial t} + Lg = 0$$
$$g(x, 0) = f(x)$$

(see Yoshida [13, Theorem 2, p. 424], for example), and defining $S(t)f := g(\cdot, t)$ we obtain a strongly continuous semigroup on $\mathcal{F} = L^2(\mathbb{R}^d)$ with infinitesimal generator $-L$.

Our first result concerns the case $s = 0$.

**Theorem 1.4.** *Let $f \in \mathcal{D}(L)$ satisfy*

$$\text{(1.12)} \qquad \int_1^\infty \|S(t)f\| \, \frac{dt}{t} < +\infty, \qquad \int_1^\infty \|S(t)Lf\| \, \frac{dt}{t} < +\infty.$$

*Then the function*

$$\text{(1.13)} \qquad u(f,y) = \int_0^\infty S(t)f \cdot e^{-y^2/4t} \, \frac{dt}{t} \qquad \text{for } y > 0$$

*satisfies*

$$-Lu(f,y) + \frac{1}{y}\frac{\partial u}{\partial y}(f,y) + \frac{\partial^2 u}{\partial y^2}(f,y) = 0 \qquad \text{for } y > 0$$

*and*

$$\lim_{y \downarrow 0} \frac{1}{\ln(y)} u(f,y) = -2f.$$

In the case $s > 0$, in analogy with the work of Stinga and Torrea [11] we may represent the solution of (1.11) via the "subordinated semigroup"

$$\text{(1.14)} \qquad u(f,y) = \frac{1}{\Gamma(s)} \int_0^\infty S(t)f \left(\frac{y^2}{4t}\right)^k e^{-y^2/4t} \, \frac{dt}{t}.$$

Indeed, if $f \in \mathcal{D}(L)$ this function satisfies (1.11) and

$$\lim_{y \downarrow 0} u(f,y) = f,$$

see Lemma 3.2. The substitution $t = y^2 q$ gives

$$\text{(1.15)} \qquad u(f,y) = \frac{1}{\Gamma(s)} \int_0^\infty S(y^2 q) f \frac{e^{-1/4q}}{(4q)^s} \, \frac{dq}{q},$$

and from this representation we again immediately see the relevance of the operator $T$, since each differentiation with respect to $\frac{1}{2}y^2$ reduces the singularity associated with the factor $q^{-s}$. This motivates the following theorem.

**Theorem 1.5.** *Let $s = k \in \mathbb{N}$, and suppose that $f, Lf, \ldots, L^k f \in \mathcal{D}(L)$ with*

$$\text{(1.16)} \qquad \int_1^\infty \|S(t)L^k f\| \, \frac{dt}{t} < +\infty.$$

*Then the function $u$ defined in (1.14) has the boundary asymptotics*

$$\text{(1.17)} \qquad \lim_{y \downarrow 0} \frac{T^k u(f,y)}{\ln(y)} = \frac{(-1)^{k-1}}{2^{k-1}\Gamma(k)} L^k f.$$

In Section 2 we derive the various constants utilized in the paper and connect them with other alternative formulations found in the literature, while in Section 3 we give proofs of the main results.

## 2. Definition of the Various Constants

At the heart of the constants in this paper are various integrations over the sphere, and therefore it is necessary to recall the Euler Beta and Gamma functions:

$$B(v,w) := \int_0^1 t^{v-1}(1-t)^{w-1} \, dt = \int_0^\infty \frac{t^{v-1}}{(t+1)^{v+w}} \, dt$$

$$\Gamma(v) := \int_0^\infty t^{v-1} e^{-t} \, dt.$$

In particular, an important identity between the Beta function and the Gamma function is the relation
$$B(v,w) = \frac{\Gamma(v)\Gamma(w)}{\Gamma(v+w)}.$$

With this notation, we can express the surface area of the $(d-1)$-dimensional sphere

$$\sigma_d = \frac{2\pi^{d/2}}{\Gamma(d/2)}. \tag{2.1}$$

We use the convention
$$\hat{f}(\xi) := \int_{\mathbb{R}^d} f(x) e^{-2\pi i \xi \cdot x}\, dx$$

for the Fourier transform. This combined with our desire that

$$((-\Delta)^s f)\hat{\,}(\xi) = (2\pi|\xi|)^{2s} \hat{f}(\xi) \tag{2.2}$$

implies that in the definition of the fractional Laplacian we take

$$c_{d,s} := \frac{2^{2s} s \Gamma(\frac{d}{2}+s)}{\pi^{\frac{d}{2}} \Gamma(1-s)}.$$

To arrive at this constant, one can use the semigroup property of the fractional Laplacian and its inverse, the Riesz potential:

$$I_\alpha f := \frac{1}{\gamma(\alpha)} \int_{\mathbb{R}^d} \frac{f(z)}{|x-z|^{d-\alpha}}\, dz,$$

where the constant (cf. Stein [10] page 117)

$$\gamma(\alpha) = \frac{\pi^{d/2} 2^\alpha \Gamma(\alpha/2)}{\Gamma(\frac{d-\alpha}{2})}.$$

Then formally one has
$$(-\Delta)^s f(x) = I_{-2s} f(x)$$

and we compute
$$\frac{1}{\gamma(-2s)} = \frac{\Gamma(\frac{d+2s}{2})}{\pi^{d/2} 2^{-2s} \Gamma(-2s/2)}$$
$$= \frac{2^{2s}\Gamma(\frac{d}{2}+s)}{\pi^{d/2}\Gamma(-s)}$$
$$= -\frac{2^{2s} s \Gamma(\frac{d}{2}+s)}{\pi^{d/2}\Gamma(1-s)},$$

where we have used $\Gamma(-s) = -\frac{\Gamma(1-s)}{s}$ for $s \in (0,1)$. From this we obtain the formal expression

$$(-\Delta)^s f(x) = \frac{2^{2s} s \Gamma(\frac{d}{2}+s)}{\pi^{d/2}\Gamma(1-s)} \int_{\mathbb{R}^d} \frac{-f(y)}{|x-y|^{d+2s}}\, dy,$$

though as the integral is not well-defined in this setting we must pursue another avenue to obtain such a representation. For example, one can obtain an admissible formula from the relation

$$(-\Delta)^s f = -\Delta(I_{2-2s}f) = I_{2-2s}(-\Delta)f.$$

Proceeding from here, an integration by parts yields

$$\begin{aligned}(-\Delta)^s f(x) &= \frac{-d+2-2s}{\gamma(2-2s)} \int_{\mathbb{R}^d} \nabla f(z) \cdot \frac{z-x}{|z-x|^{d+2s}} \, dz \\ &= \frac{-d+2-2s}{\gamma(2-2s)} \int_{\mathbb{R}^d} \nabla [f(z)-f(x)] \cdot \frac{z-x}{|z-x|^{d+2s}} \, dz \\ &= -\frac{-d+2-2s}{\gamma(2-2s)} \, p.v. \int_{\mathbb{R}^d} [f(z)-f(x)] \left[ \frac{d}{|z-x|^{d+2s}} + \frac{-d-2s}{|z-x|^{d+2s}} \right] dz \\ &= 2s\frac{d-2+2s}{\gamma(2-2s)} \, p.v. \int_{\mathbb{R}^d} \frac{f(x)-f(z)}{|x-z|^{d+2s}} \, dz.\end{aligned}$$

The definition of $\gamma$ and the factorial extension property of the Gamma function $v\Gamma(v) = \Gamma(v+1)$ then yields

$$\begin{aligned}2s\frac{d-2+2s}{\gamma(2-2s)} &= 2^{2s}\frac{(d-2+2s)}{2} \frac{\Gamma(\frac{d-(2-2s)}{2})}{\pi^{d/2}2^{2-2s}\Gamma((2-2s)/2)} \\ &= \frac{2^{2s}s\Gamma(\frac{d-(2-2s)}{2}+1)}{\pi^{d/2}\Gamma(1-s)},\end{aligned}$$

which simplifies to give $c_{d,s}$ as above.

The Poisson-Bessel kernel of the equation (1.2) has been computed by Marias [6, p. 320, equation (3.1)] or Caffarelli and Silvestre [1, Section 2]. In particular we can write the solution $u$ via the formulae

$$\begin{aligned}u(x,y) &= \frac{\Gamma(\frac{d}{2}+s)}{\pi^{d/2}\Gamma(s)} \, y^{2s} \int_{\mathbb{R}^d} \frac{f(x-z)}{(|z|^2+y^2)^{d/2+s}} \, dz \\ &= \frac{\Gamma(\frac{d}{2}+s)}{\pi^{d/2}\Gamma(s)} \, y^{2s} \int_{\mathbb{R}^d} \frac{f(z)}{(|x-z|^2+y^2)^{d/2+s}} \, dz.\end{aligned}$$

To understand the value of the constant in the normalization, one use the change of variables $h = -z/y$ in the first formula for $u$ to find

$$u(x,y) = \frac{\Gamma(\frac{d}{2}+s)}{\pi^{d/2}\Gamma(s)} \int_{\mathbb{R}^d} \frac{f(x+yh)}{(|h|^2+1)^{d/2+s}} \, dh.$$

Thus, in order to show that this integral is an approximation of the identity for any $s > 0$, we must verify that

$$(2.3) \qquad \frac{\Gamma(\frac{d}{2}+s)}{\pi^{d/2}\Gamma(s)} = \left( \int_{\mathbb{R}^d} \frac{1}{(|h|^2+1)^{d/2+s}} \, dh \right)^{-1}.$$

However, integrating in spherical coordinates and changing variables $w = t^2$ we find

$$\begin{aligned}\int_{\mathbb{R}^d} \frac{1}{(|h|^2+1)^{d/2+s}} \, dh &= \sigma_d \int_0^\infty \frac{t^{d-1}}{(t^2+1)^{d/2+s}} \, dt \\ &= \frac{\sigma_d}{2} \int_0^\infty \frac{w^{d/2-1}}{(w+1)^{d/2+s}} \, dw \\ &= \frac{\sigma_d}{2} B\left(\frac{d}{2},s\right),\end{aligned}$$

where $B$ is the beta function as defined earlier in this section. Using its relation with the Gamma function and the formula for $\sigma_d$ we find

$$\frac{\sigma_d}{2} B\left(\frac{d}{2},s\right) = \frac{\Gamma(\frac{d}{2}+s)}{\pi^{d/2}\Gamma(s)},$$

which is precisely (2.3).

2.1. **Computation of the constant of Caffarelli and Silvestre.** We can now compute explicitly the constant arising when we evaluate the limits (1.3) and (1.7).

Let us first consider the case $s \in (0,1)$, where we will evaluate the limit

$$\tilde{c}_s(-\Delta)^s f(x) = \lim_{y \downarrow 0} \frac{u(x,y) - u(x,0)}{y^{2s}/2s}.$$

We use the second formula in the above definition of $u$ and evaluate the right hand side of the formula (1.3) to find

$$\lim_{y \downarrow 0} \frac{u(x,y) - u(x,0)}{y^{2s}/2s} = \lim_{y \downarrow 0} 2s \frac{\Gamma(\frac{d}{2} + s)}{\pi^{d/2}\Gamma(s)} \int_{\mathbb{R}^d} \frac{f(z) - f(x)}{(|x-z|^2 + y^2)^{d/2+s}} \, dz$$

$$= \lim_{y \downarrow 0} -2s \frac{\Gamma(\frac{d}{2} + s)}{\pi^{d/2}\Gamma(s)} \int_{\mathbb{R}^d} \frac{f(x) - f(z)}{(|x-z|^2 + y^2)^{d/2+s}} \, dz$$

$$= \frac{-2s}{c_{d,s}} \frac{\Gamma(\frac{d}{2} + s)}{\pi^{d/2}\Gamma(s)} (-\Delta)^s f(x),$$

which leads to the constant

$$\tilde{c}_s = -2s \times \frac{\Gamma(\frac{d}{2} + s)}{\pi^{d/2}\Gamma(s)} \times \frac{\pi^{\frac{d}{2}}\Gamma(1-s)}{2^{2s}s\Gamma(\frac{d}{2}+s)}$$

$$= -\frac{\Gamma(1-s)}{2^{2s-1}\Gamma(s)},$$

as in (1.3).

For the case $s \geq 1$, we proceed as follows. First, from the proof of Theorem 1.2 we find a formula for repeated application of the operator $T$ in equation 3.2. In particular, if $u$ is the solution of (1.2) for $s \in \mathbb{R}^+$, then for $k = \lfloor s \rfloor$, $\sigma = s - k$ we have

$$T^k u(x,y) = \frac{(-1)^k \Gamma(\sigma)}{2^k \Gamma(s)} C_{d,\sigma} \int_{\mathbb{R}^d} \frac{(-\Delta)^k f(x-z)}{(|z|^2 + y^2)^{d/2+\sigma}} \, dz.$$

Therefore combined with the computation for $\sigma \in (0,1)$ we find

$$\lim_{y \downarrow 0} 2\sigma \frac{T^k u(x,y) - T^k u(x,0)}{y^{2\sigma}} = \frac{(-1)^k \Gamma(\sigma)}{2^k \Gamma(s)} \tilde{c}_\sigma (-\Delta)^s f(x).$$

In particular, for all $s \in \mathbb{R}_+ \setminus \mathbb{N}$ one has

(2.4) $$\tilde{c}_s = \frac{(-1)^k \Gamma(\sigma)}{2^k \Gamma(s)} \times \left[ -\frac{\Gamma(1-\sigma)}{2^{2\sigma-1}\Gamma(\sigma)} \right] = \frac{(-1)^{k-1} 2^k \Gamma(1-\sigma)}{2^{2s-1}\Gamma(s)}.$$

Observe that the formula (2.4) coincides with the constant appearing in Theorem 1.2 when $s = k \in \mathbb{N}$ (which is the case $\sigma = 0$). This will be established in the proof of Theorem 1.2 in Section 3.

2.2. **Comparison with the constant of Chang and Gonzalez.** In the preceding section we find one value for the constant associated to the boundary limit (1.7). Meanwhile, in the work of Chang and Gonzalez [2, p. 1420] they have shown that for $s \in (0, d/2)$ one has

(2.5) $$(-\Delta)^s f(x) = \frac{2^{2s}}{2\sigma} \frac{\Gamma(s)}{\Gamma(-s)} \frac{1}{2^k s(s-1)(s-2)\cdots(s-k+1))}$$

$$\cdot \lim_{y \downarrow 0} y^{1-2\sigma} \frac{\partial}{\partial y} T^k(u^s)(x,y).$$

Here, as before $k = \lfloor s \rfloor$ and $\sigma = s - k \in (0,1)$. Let us note that on p. 1420, the incorrect constant for $A_m$ is given, and one should look on its definition on p. 1419 to see the missing factor of $s$ that we have taken in the above formula.

Let us here reconcile the two formula and their associated constants. To this end, first let us observe that the boundary operator

$$\lim_{y \downarrow 0} y^{1-2\sigma} \frac{\partial u}{\partial y}(x,y) = \lim_{y \downarrow 0} \frac{\partial u}{\partial (y^{2\sigma}/2\sigma)}(x,y).$$

is equivalent to evaluating the limit

$$\lim_{y \downarrow 0} \frac{u(x,y) - u(x,0)}{y^{2\sigma}/2\sigma},$$

provided the limits exist. In particular, in our setting they do exist and thus the two boundary operators are the same. As a consequence one can alternatively state our result as

$$(2.6) \qquad (-\Delta)^s f(x) = \frac{(-1)^{k-1} 2^{2s-1} \Gamma(s)}{2^k \Gamma(1-\sigma)} \lim_{y \downarrow 0} y^{1-2\sigma} \frac{\partial}{\partial y} T^k(u^s)(x,y).$$

It therefore suffices to show that the constants in equations (2.5) and (2.6) are the same. This will be established once we have checked that

$$(2.7) \qquad \frac{2^{2s}}{2\sigma} \frac{\Gamma(s)}{\Gamma(-s)} \frac{1}{2^k s(s-1)(s-2)\cdots(s-k+1)} \times \left[ -\frac{(-1)^k 2^k \Gamma(1-\sigma)}{2^{2s-1} \Gamma(s)} \right] = 1.$$

But the factor of $(-1)^k$, the falling factorial in the denominator, and $\Gamma(-s)$ can be combined as

$$\frac{(-1)^k}{\Gamma(-s)} \frac{1}{s(s-1)(s-2)\cdots(s-k+1))} = \frac{1}{\Gamma(-s+k)} = \frac{1}{\Gamma(-\sigma)}.$$

Inserting this into (2.7), the factors of $\Gamma(s)$ and powers of 2 in the numerator and denominator cancel to leave

$$\frac{\Gamma(1-\sigma)}{-\sigma \Gamma(-\sigma)} = 1,$$

as required, where we have again used the identity $v\Gamma(v) = \Gamma(v+1)$.

## 3. Proofs of the Main Results

*Proof of Theorem 1.1.* Let $f \in C_c^\infty(\mathbb{R}^d)$ and define the function $u$ by

$$u(x,y) = \frac{\Gamma(d/2)}{\pi^{d/2}} \int_{\mathbb{R}^d} \frac{f(x-z)}{(|z|^2 + y^2)^{d/2}} \, dz.$$

Then if we define

$$K(x,y) := \frac{1}{(|x|^2 + y^2)^{d/2}},$$

we have that $u(x,y) = \Gamma(d/2)\pi^{-d/2}(K(\cdot, y) * f)(x)$, and so it suffices to check that $K$ satisfies

$$\Delta_{x,y} K(x,y) + \frac{1}{y} \frac{\partial K}{\partial y}(x,y) = 0.$$

But we have

$$\nabla_x K(x,y) = -d \frac{x}{(|x|^2 + y^2)^{d/2+1}}$$

$$\frac{\partial K}{\partial y}(x,y) = -d \frac{y}{(|x|^2 + y^2)^{d/2+1}},$$

so that

$$\Delta_{x,y}K(x,y) = \frac{-d(d+1)}{(|x|^2+y^2)^{d/2+1}} + \frac{-d(-d-2)(|x|^2+y^2)}{(|x|^2+y^2)^{d/2+2}}$$
$$= \frac{-d(d+1) - d(-d-2)}{(|x|^2+y^2)^{d/2+1}}$$
$$= \frac{d}{(|x|^2+y^2)^{d/2+1}},$$

which implies the claim.

For the boundary limit, we make the change of variables $h = \frac{-z}{y}$ to obtain

$$u(x,y) = \frac{\Gamma(d/2)}{\pi^{d/2}} \int_{\mathbb{R}^d} \frac{f(x+yh)}{(|h|^2+1)^{d/2}}\, dh.$$

Thus, if $\operatorname{supp} f \subset B(0,R)$, then we have

$$\lim_{y \downarrow 0} \frac{-u(x,y)}{\ln(y)} = \lim_{y \downarrow 0} \frac{-1}{\ln(y)} \int_{B(x,R/y)} \frac{f(x+yh)}{(|h|^2+1)^{d/2}}\, dh.$$

Now notice that the constant $\Gamma(d/2)/\pi^{d/2}$ coincides with $2/\sigma_d$. Thus if we can show that

(3.1) $$\lim_{y \downarrow 0} \frac{-1}{\ln(y)} \frac{1}{\sigma_d} \int_{B(x,R/y)} \frac{1}{(|h|^2+1)^{d/2}}\, dh = 1,$$

then we will have

$$\lim_{y \downarrow 0} \frac{-u(x,y)}{\ln(y)} = \lim_{y \downarrow 0} \frac{-1}{\ln(y)} \frac{2}{\sigma_d} \int_{B(x,R/y)} \frac{f(x+yh) - f(x)}{(|h|^2+1)^{d/2}}\, dh + 2f(x),$$

and since

$$\left| \frac{1}{\ln(y)} \int_{B(x,R/y)} \frac{f(x+yh) - f(x)}{(|h|^2+1)^{d/2}}\, dh \right| \leq \left| \frac{1}{\ln(y)} \right| \|\nabla f\|_\infty \int_{B(x,R/y)} \frac{y|h|}{(|h|^2+1)^{d/2}}\, dh$$
$$= \left| \frac{1}{\ln(y)} \right| \|\nabla f\|_\infty \int_{B(x,R)} \frac{|z|}{(|z|^2+y^2)^{d/2}}\, dz$$

tends to zero as $y \downarrow 0$, the result will be demonstrated.

Thus, it remains to show (3.1). But

$$\int_{B(x,R/y)} \frac{1}{(|h|^2+1)^{d/2}}\, dh = \int_0^{R/y} \int_{S^{d-1}} \frac{1}{(|x-\rho\theta|^2+1)^{d/2}} \rho^{d-1}\, d\mathcal{H}^{d-1}(\theta)\, d\rho,$$

so that

$$\lim_{y \downarrow 0} \frac{-1}{\ln(y)} \int_0^{R/y} \int_{S^{d-1}} \frac{\rho^{d-1}}{(|x-\rho\theta|^2+1)^{d/2}}\, d\mathcal{H}^{d-1}(\theta)\, d\rho$$

is an indefinite form to which we can apply L'Hopital's rule. Thus,

$$\lim_{y \downarrow 0} \frac{-1}{\ln(y)} \int_0^{R/y} \int_{S^{d-1}} \frac{\rho^{d-1}}{(|x-\rho\theta|^2+1)^{d/2}}\, d\mathcal{H}^{d-1}(\theta)\, d\rho$$
$$= \lim_{y \downarrow 0} y \frac{R}{y^2} \int_{S^{d-1}} \frac{(R/y)^{d-1}}{(|x-(R/y)\theta|^2+1)^{d/2}}\, d\mathcal{H}^{d-1}(\theta)$$
$$= \lim_{y \downarrow 0} \int_{S^{d-1}} \frac{1}{(|(y/R)x - \theta|^2 + (y/R)^2)^{d/2}}\, d\mathcal{H}^{d-1}(\theta)$$
$$= \sigma_d \qquad \qquad \square$$

*Proof of Theorem 1.2.* Let $s = k \in \mathbb{N}$ and $f \in C_c^\infty(\mathbb{R}^d)$. In this regime of $s$, the boundary value problem (1.2) is well-posed and the solution is given by

$$u(x,y) = \frac{\Gamma(\frac{d}{2}+k)}{\pi^{d/2}\Gamma(k)} y^{2k} \int_{\mathbb{R}^d} \frac{f(x-z)}{(|z|^2+y^2)^{d/2+k}} dz.$$

Then we will compute $T^k u$. For convenience of reference let us perform the computation for any $s \geq 1$, where solution of (1.2) is given by

$$u(x,y) = \frac{\Gamma(\frac{d}{2}+s)}{\pi^{d/2}\Gamma(s)} \int_{\mathbb{R}^d} \frac{f(x+yh)}{(|h|^2+1)^{d/2+s}} dh.$$

Then Lebesgue's dominated convergence theorem justifies the differentiation

$$\frac{\partial u}{\partial y}(x,y) = \frac{\Gamma(\frac{d}{2}+s)}{\pi^{d/2}\Gamma(s)} \int_{\mathbb{R}^d} \frac{\nabla f(x+yh) \cdot h}{(|h|^2+1)^{d/2+s}} dh,$$

and therefore

$$Tu(x,y) = \frac{\Gamma(\frac{d}{2}+s)}{\pi^{d/2}\Gamma(s)} \frac{1}{y} \int_{\mathbb{R}^d} \frac{\nabla f(x+yh) \cdot h}{(|h|^2+1)^{d/2+s}} dh.$$

Now,

$$\frac{h}{(|h|^2+1)^{d/2+s}} = \frac{-1}{2} \frac{1}{d/2+s-1} \nabla \frac{1}{(|h|^2+1)^{d/2+s-1}},$$

and thus performing an integration by parts we find

$$Tu(x,y) = \frac{\Gamma(\frac{d}{2}+s)}{\pi^{d/2}\Gamma(s)} \frac{-1}{2} \frac{1}{d/2+s-1} \int_{\mathbb{R}^d} \frac{-\Delta f(x+yh)}{(|h|^2+1)^{d/2+s-1}} dh.$$

Write $s = k + \sigma$, where $k = \lfloor s \rfloor \in \mathbb{N}$ and $\sigma \in [0,1)$. We can continue this process $k$ times to find

$$T^k u(x,y) = \frac{\Gamma(\frac{d}{2}+s)}{\pi^{d/2}\Gamma(s)} \left(\frac{-1}{2}\right)^k \frac{1}{d/2+s-1} \cdots \frac{1}{d/2+s-k}$$
$$\cdot \int_{\mathbb{R}^d} \frac{(-\Delta)^k f(x-z)}{(|z|^2+y^2)^{d/2+s-k}} dz.$$

But now using $v\Gamma(v) = \Gamma(v+1)$ and $s - k = \sigma$ we have

$$\Gamma(\tfrac{d}{2}+s) \times \frac{1}{d/2+s-1} \cdots \frac{1}{d/2+s-k} = \Gamma(\tfrac{d}{2}+\sigma).$$

Thus we obtain

(3.2) $$T^k u(x,y) = \frac{(-1)^k}{2^k \Gamma(s)} \frac{\Gamma(\frac{d}{2}+\sigma)}{\pi^{d/2}} \int_{\mathbb{R}^d} \frac{(-\Delta)^k f(x-z)}{(|z|^2+y^2)^{d/2+\sigma}} dz.$$

With these preparations, we can now take $s = k$, $\sigma = 0$ and apply Theorem 1.1 to deduce that

$$\lim_{y \downarrow 0} \frac{1}{\ln(y)} T^k u(x,y) = \frac{(-1)^k}{2^k \Gamma(k)} \lim_{y \downarrow 0} \frac{1}{\ln(y)} \frac{\Gamma(\frac{d}{2})}{\pi^{d/2}} \int_{\mathbb{R}^d} \frac{(-\Delta)^k f(x-z)}{(|z|^2+y^2)^{d/2}} dz$$
$$= \frac{(-1)^k}{2^k \Gamma(k)} \times \left(-2(-\Delta)^k f(x)\right). \qquad \square$$

*Proof of Proposition 1.3.* We will show that

(3.3) $$\frac{\Gamma(d/2)}{\pi^{d/2}} \int_{\mathbb{R}^d} \frac{f(z)}{(|x-z|^2+y^2)^{d/2}} dz = \int_0^\infty e^{t\Delta} f(x) e^{-\frac{y^2}{4t}} \frac{dt}{t}.$$

Let us consider the integral on the right hand side of (3.3). Using the fact that the heat semigroup is given by integration against the heat kernel, we have

$$e^{t\Delta} f(x) = \frac{1}{(4\pi t)^{d/2}} \int_{\mathbb{R}^d} f(z) e^{-\frac{|x-z|^2}{4t}} \, dz.$$

This expression and Fubini's theorem imply that

$$(3.4) \qquad \int_0^\infty e^{t\Delta} f(x) e^{-\frac{y^2}{4t}} \frac{dt}{t} = \int_{\mathbb{R}^d} f(z) \int_0^\infty \frac{1}{(4\pi t)^{d/2}} e^{-\frac{|x-z|^2}{4t}} e^{-\frac{y^2}{4t}} \frac{dt}{t} \, dz.$$

Considering the inner integral and changing variables $v = -(|x-z|^2 + y^2)/4t$, we find

$$\int_0^\infty \frac{1}{(4\pi t)^{d/2}} e^{-\frac{|x-z|^2}{4t}} e^{-\frac{y^2}{4t}} \frac{dt}{t} = \frac{1}{(|x-z|^2 + y^2)^{d/2}} \frac{1}{\pi^{d/2}} \int_0^\infty s^{d/2} e^{-s} \frac{ds}{s}$$
$$= \frac{1}{(|x-z|^2 + y^2)^{d/2}} \frac{\Gamma(d/2)}{\pi^{d/2}},$$

so combining with (2.1) and (3.4) recovers the left hand side of equation (3.3). □

Theorem 1.4 asserts that the function $u(f, y)$ defined by (1.13) solves a PDE and has logarithmic boundary behavior. We prove the first assertion in a separate lemma. The operator $L$ will arise through differentiation with respect to $t$ in the integral defining $u(f, y)$, using the identity

$$(3.5) \qquad \frac{\partial}{\partial t} S(t) f = -L S(t) f = -S(t) L f,$$

see for instance [8, Lemma 12.11].

**Lemma 3.1.** *Under the hypotheses of Theorem 1.4, $u$ satisfies*

$$-L u(f, y) + \frac{1}{y} \frac{\partial u}{\partial y}(f, y) + \frac{\partial^2 u}{\partial y^2}(f, y) = 0 \qquad \text{for } y > 0.$$

The proof is based on the observation that the kernel

$$K(t, y) = \frac{1}{t} e^{-y^2/4t}$$

satisfies

$$(3.6) \qquad \left( \frac{1}{y} \frac{\partial}{\partial y} + \frac{\partial^2}{\partial y^2} \right) K(t, y) = \frac{\partial}{\partial t} K(t, y) \qquad \text{for } t > 0, y > 0,$$

together with an integration by parts and the identity (3.5). Formalizing this intuition requires appropriate care for the interchange of integration and differentiation.

*Proof.* It is readily verified that $K(t, y)$, $\frac{\partial K}{\partial t}(t, y)$ and $\frac{\partial^2 K}{\partial y^2}$ can all be bounded, uniformly over $y$ in a compact subset of $(0, \infty)$, by a constant times $\min\{1/t, 1\}$. Together with the boundedness of $\|S(t)f\|$ and the assumption (1.12), it follows from the dominated convergence theorem for Bochner integrals that

$$\frac{\partial^m}{\partial y^m} u(f, y) = \int_0^\infty S(t) f \cdot \frac{\partial^m K}{\partial y^m}(t, y) \, dt \qquad m = 1, 2,$$

and therefore

$$\left( \frac{1}{y} \frac{\partial}{\partial y} + \frac{\partial^2}{\partial y^2} \right) u(f, y) = \int_0^\infty S(t) f \left( \frac{1}{y} \frac{\partial}{\partial y} + \frac{\partial^2}{\partial y^2} \right) K(t, y) \, dt.$$

Using (3.6),

$$\left( \frac{1}{y} \frac{\partial}{\partial y} + \frac{\partial^2}{\partial y^2} \right) u(f, y) = \int_0^\infty S(t) f \frac{\partial K}{\partial t}(t, y) \, dt.$$

For any finite $t_0 < \infty$, we can integrate by parts to find

$$\text{(3.7)} \quad \int_0^{t_0} S(t)f \, \frac{\partial K}{\partial t}(t,y) \, dt = S(t)f \cdot K(t,y)\big|_{t=0}^{t=t_0} - \int_0^{t_0} \left(\frac{\partial}{\partial t} S(t)f\right) K(t,y) \, dt$$

$$= S(t_0)f \cdot K(t_0, y) + \int_0^{t_0} LS(t)f \cdot K(t,y) \, dt$$

by (3.5). As $t_0 \to \infty$, the first term on the right-hand side of (3.7) vanishes by (1.10), and the integral converges by (1.12). The integral on the left-hand side of (3.7) converges by (1.10) and the absolute integrability of $\int_0^\infty \frac{\partial K}{\partial t}(t,y) \, dt$. We conclude that

$$\left(\frac{1}{y}\frac{\partial}{\partial y} + \frac{\partial^2}{\partial y^2}\right) u(f, y) = \int_0^\infty LS(t)f \cdot K(t,y) \, dt.$$

We wish to recognize the last integral as $Lu(f, y)$. To prove this we use the definition (1.9) of $Lu(f, y)$ as the limit of a difference quotient:

$$Lu(f, y) = \lim_{h \downarrow 0} \frac{1}{h} \left(S(h)u(f,y) - u(f,y)\right)$$

$$= \lim_{h \downarrow 0} \frac{1}{h} \int_0^\infty K(t,y) \left(S(h)S(t)f - S(t)f\right) dt$$

$$= \lim_{h \downarrow 0} \frac{1}{h} \int_0^\infty K(t,y) \left(S(t+h)f - S(t)f\right) dt$$

$$= \lim_{h \downarrow 0} \frac{1}{h} \int_0^\infty K(t,y) \, dt \int_t^{t+h} LS(\tau)f \, d\tau$$

$$= \lim_{h \downarrow 0} \frac{1}{h} \int_0^\infty \int_0^\infty \mathbb{1}_{\{t \leq \tau\}} \mathbb{1}_{\{t+h \geq \tau\}} LS(\tau)f \cdot K(t,y) \, dt \, d\tau,$$

where we have used the fact that $S(h)$ is a bounded linear operator on $\mathcal{F}$ to bring $S(h)$ inside the integral, and the identity (3.5) to write the difference quotient in terms of the derivative.

For each fixed $y > 0$, $K(t,y) = O(1)$ for $t \leq 1$ and $K(t,y) = O(t^{-1})$ for $t \geq 1$. It follows that

$$\text{(3.8)} \quad \int_{\max\{\tau - h, 0\}}^{\tau} \frac{1}{h} K(t,y) \, dt \leq C(y) \min\{1, \tau^{-1}\}$$

for some $C(y) < +\infty$, uniformly over $h \in (0, 1)$. By (1.12),

$$\int_0^\infty \|LS(\tau)f\| \int_{\max\{\tau-h, 0\}}^{\tau} K(t,y) \, dt \, d\tau < +\infty.$$

Thus, Fubini's theorem yields

$$Lu(f, y) = \lim_{h \downarrow 0} \int_0^\infty LS(\tau)f \int_{\max\{\tau-h, 0\}}^{\tau} \frac{1}{h} K(t,y) \, dt \, d\tau.$$

But now the inner integral converges to $K(\tau, y)$ pointwise for each $\tau > 0$, while the bounds (1.10), (1.12) and (3.8) allow us to invoke the dominated convergence for Bochner integrals to conclude the desired result. □

We next prove the claim in the introduction that a similar result holds for the case $s > 0$.

**Lemma 3.2.** *Suppose $f \in \mathcal{D}(L)$. Then the function $u$ defined in (1.14) satisfies*

$$\text{(3.9)} \quad -Lu(f,y) + \frac{1-2k}{y} \frac{\partial u}{\partial y}(f,y) + \frac{\partial^2 u}{\partial y^2}(f,y) = 0 \qquad \text{for } y > 0$$

and

(3.10) $$\lim_{y \downarrow 0} u(f, y) = f.$$

*Proof.* To show that $u$ satisfies (3.9), note that the kernel

$$K(t, y) = \frac{1}{t\Gamma(k)} \left(\frac{y^2}{4t}\right)^k e^{-y^2/4t}$$

satisfies $(\frac{1-2k}{y}\frac{\partial}{\partial y} + \frac{\partial^2}{\partial y^2})K(t, y) = \frac{\partial}{\partial t}K(t, y)$ for $y > 0$. The kernel $K$ is better behaved than in Lemma 3.1: the integrals $\int_0^\infty K(t, y)\, dt$, $\int_0^\infty \frac{\partial K}{\partial y}(t, y)\, dt$ and $\int_0^\infty \frac{\partial^2 K}{\partial y^2}(t, y)\, dt$ all converge absolutely since $k > 0$. We can therefore repeat the argument from the preceding proof, using only the boundedness of $\|S(t)f\|$ and $\|S(t)Lf\|$ from (1.10).

To show (3.10), substitute $t = y^2 q$ in the definition (1.14) of $u(f, y)$ to obtain

(3.11) $$u(f, y) = \frac{1}{\Gamma(k)} \int_0^\infty S(y^2 q)f\, \frac{e^{-1/4q}}{(4q)^k}\, \frac{dq}{q}.$$

Noting that

$$\int_0^\infty (4q)^{-k} e^{-1/4q} \frac{dq}{q} = \Gamma(k) \qquad \text{for all } k > 0,$$

that the strong continuity of $S(t)$ implies that

$$\lim_{y \downarrow 0} S(y^2 q)f = f$$

for any $q \geq 0$ and any $f$, and recalling the boundedness from (1.10), we can apply the dominated convergence theorem for Bochner integrals to conclude (3.10). □

Finally, it will be convenient to treat Theorems 1.4 and 1.5 in parallel, and therefore we first handle the somewhat delicate analysis of repeated application of the operator $T$.

**Lemma 3.3.** *Under the hypotheses of Theorem 1.5,*

(3.12) $$T^k u(f, y) = \frac{(-1)^k}{2^k \Gamma(k)} \int_0^\infty S(y^2 q) L^k f \cdot e^{-1/4q} \frac{dq}{q}.$$

*Proof.* Consider first $k > 1$. Recalling that $T$ is a factor of two times differentiation with respect to $y^2$, we can use (3.11) to write

$$Tu(f, y) = \lim_{h \downarrow 0} \frac{2}{\Gamma(k)} \int_0^\infty \frac{S((y^2 + h)q)f - S(y^2 q)f}{h} \cdot \frac{e^{-1/4q}}{(4q)^k} \frac{dq}{q}$$

(3.13) $$= \lim_{h \downarrow 0} \frac{2}{4^k \Gamma(k)} \int_0^\infty \frac{S(y^2 q + hq)f - S(y^2 q)f}{hq} \cdot q \frac{e^{-1/4q}}{q^k} \frac{dq}{q}.$$

As $h \downarrow 0$, the integrand converges pointwise to $-S(y^2 q)Lf \cdot e^{-1/4q} q^{-k}$ by (3.5). Using (1.10), we can apply the dominated convergence theorem for Bochner integrals to obtain

$$Tu(f, y) = -\frac{2}{4^k \Gamma(k)} \int_0^\infty S(y^2 q)Lf \cdot \frac{e^{-1/4q}}{q^{k-1}} \frac{dq}{q}.$$

Because $f, \ldots, L^{k-2}f \in \mathcal{D}(L)$, we can repeat this argument as long as the limiting kernel remains integrable, a total of $k - 1$ times, to obtain

(3.14) $$T^{k-1} u(f, y) = \frac{(-1)^{k-1}}{2^{k+1} \Gamma(k)} \int_0^\infty S(y^2 q) L^{k-1} f \cdot \frac{e^{-1/4q}}{q} \frac{dq}{q}.$$

For the final application of $T$, we argue as in (3.13) and the proof of Lemma 3.1:

$$T^k u(f,y) = \lim_{h\downarrow 0} \frac{2(-1)^{k-1}}{2^{k+1}\Gamma(k)} \int_0^\infty \frac{S(y^2 q + hq)L^{k-1}f - S(y^2 q)L^{k-1}f}{hq} \cdot e^{-1/4q} \frac{dq}{q}$$

$$= \lim_{h\downarrow 0} \frac{(-1)^{k-1}}{2^k \Gamma(k)} \int_0^\infty \frac{1}{hq} \int_{y^2 q}^{y^2 q + hq} (-S(t)L^k f) \cdot e^{-1/4q} \, dt \, \frac{dq}{q}$$

(3.15)
$$= \lim_{h\downarrow 0} \frac{(-1)^k}{2^k \Gamma(k)} \int_0^\infty S(t) L^k f \cdot \frac{1}{h} \int_{t/(y^2+h)}^{t/y^2} e^{-1/4q} \frac{dq}{q^2} \, dt.$$

The mean value theorem for integrals gives

$$\frac{1}{h} \int_{t/(y^2+h)}^{t/y^2} e^{-1/4q} \frac{dq}{q^2} = \frac{1}{h} e^{-1/4q^*} \int_{t/(y^2+h)}^{t/y^2} \frac{dq}{q^2} = \frac{e^{-1/4q^*}}{t}$$

for some $q^* = q^*(t,h,y) \in (\frac{t}{y^2+h}, \frac{t}{y^2})$, and so

(3.16)
$$T^k u(f,y) = \lim_{h\downarrow 0} \frac{(-1)^k}{2^k \Gamma(k)} \int_0^\infty S(t) L^k f \cdot \frac{e^{-1/4q^*}}{t} \, dt.$$

But now $q^* \to t/y^2$ as $h \downarrow 0$ so that the integrand in (3.16) converges pointwise to $S(t) L^k f \cdot t^{-1} e^{-y^2/4t}$. Meanwhile

$$\frac{e^{-1/4q^*}}{t} \leq \frac{1}{t} e^{-y^2/4t} \leq C(y) \min\{1, 1/t\}$$

so that (1.10) and (1.16) justifies the application of the dominated convergence theorem for Bochner integrals and the proof is complete. □

*Proof of Theorems 1.4 and 1.5.* To treat the cases $k = 0$ and $k \in \mathbb{N}$ in parallel, set $\tilde{\Gamma}(k) = \Gamma(k)$ if $k > 0$ and $\tilde{\Gamma}(0) = 1$, so that

$$u(f,y) = \frac{1}{\tilde{\Gamma}(k)} \int_0^\infty S(t) f \left(\frac{y^2}{4t}\right)^k e^{-y^2/4t} \frac{dt}{t}$$

in both cases. We begin with the formula established in Lemma 3.3 for $k \in \mathbb{N}$:

(3.17)
$$T^k u(f,y) = \frac{(-1)^k}{2^k \tilde{\Gamma}(k)} \int_0^\infty S(y^2 q) f \cdot e^{-1/4q} \frac{dq}{q}.$$

We note that (3.17) is the definition of $u$ in the case $k = 0$.

Since the integral $\int_0^\infty e^{-1/4q} dq/q$ is divergent, we cannot set $y = 0$. Instead, write

$$S(t) L^k f = e^{-t} L^k f + \left[ S(t) L^k f - L^k f + (1 - e^{-t}) L^k f \right].$$

Since $L^k f \in \mathcal{D}(L)$, the function $t \mapsto S(t) L^k f$ is differentiable at $t = 0$. In particular, $\|S(t) L^k f - L^k f\| = O(t)$ for $t \leq 1$, while (1.10) gives $\|S(t) L^k f - L^k f\| = O(1)$ for $t \geq 1$. These two asymptotic statements can be combined to say that $\|S(t) L^k f - L^k f\| = O(1 - e^{-t})$ for all $t \geq 0$, and hence

$$S(t) L^k f = e^{-t} L^k f + O(1 - e^{-t}),$$

where the $O(\cdot)$ term means a function $g(t) : [0, \infty) \to \mathcal{F}$ with $\|g(t)\| = O(1 - e^{-t})$. Inserting into (3.12) and substituting back in terms of $t$,

$$T^k u(f,y) = L^k f \cdot \frac{(-1)^k}{2^k \tilde{\Gamma}(k)} \int_0^\infty e^{-t} e^{-y^2/4t} \frac{dt}{t} + O(1) \int_0^\infty (1 - e^{-t}) e^{-y^2/4t} \frac{dt}{t}.$$

The second integral converges even when $y = 0$, while the first reduces to a modified Bessel function:

$$T^k u(f,y) = L^k f \cdot \frac{(-1)^k}{2^{k-1} \tilde{\Gamma}(k)} K_0(y) + O(1)$$

by [4, 8.432.6]. Applying the asymptotics $K_0(y) = -\ln(y) + O(y\ln(y))$ from [4, 8.447.3], we conclude that

$$T^k u(f, y) = \frac{(-1)^{k-1}}{2^{k-1}\tilde{\Gamma}(k)} L^k f \cdot \ln(y) + O(1/\ln(y)) + O(y)$$

as $y \downarrow 0$, and thus the result is demonstrated. $\square$

## Acknowledgements

The authors would like to thank Cindy Chen for her help in obtaining some of the articles utilized in this research and the Technion for the stimulating environment that led to the initiation of this collaboration. D.S. would like to thank the University of Auckland for its warm hospitality during which some of this research was undertaken. D.S. is supported in part by the Taiwan Ministry of Science and Technology under research grants 103-2115-M-009-016-MY2 and 105-2115-M-009-004-MY2. J.G. would like to thank the National Chiao Tung University for its warm hospitality during which some of this research was undertaken. J.G. is supported in part by the Marsden Fund Council from New Zealand Government funding, managed by the Royal Society of New Zealand.

Department of Statistics
University of Auckland
Private Bag 92019, Victoria Street West
Auckland 1142, New Zealand
*E-mail address*: jesse.goodman@auckland.ac.nz

Department of Applied Mathematics
National Chiao Tung University
1001 Ta Hsueh Road.
Hsinchu, Taiwan R.O.C.

National Center for Theoretical Sciences
National Taiwan University
No. 1 Sec. 4 Roosevelt Rd.
Taipei, 106, Taiwan R.O.C.
*E-mail address*: `dspector@math.nctu.edu.tw`